\documentclass{tran-l}
\numberwithin{equation}{section}
\newtheorem{thm}{Theorem}[section]
\newtheorem{lem}[thm]{Lemma}
\newtheorem{prop}[thm]{Proposition}
\newtheorem{cor}[thm]{Corollary}

\theoremstyle{definition}
\theoremstyle{remark}
\newtheorem*{rem}{Remark}
\newtheorem*{concl}{Conclusion}

\newcommand{\W}{\mathcal{W}}
\newcommand{\K}{\mathcal{K}}
\newcommand{\AK}{\mathcal{AK}}

\newcommand{\thmref}[1]{Theorem~\ref{#1}}
\newcommand{\secref}[1]{\S\ref{#1}}
\newcommand{\lemref}[1]{Lemma~\ref{#1}}
\begin{document}

\title[Tangent Bundles with Sasaki Metric and Almost $(H,G)$-Structure]
{Tangent Bundles with Sasaki Metric and Almost Hypercomplex Pseudo-Hermitian
Structure}

\author{Mancho Manev}

\address{University of Plovdiv, Faculty of Mathematics and
Informatics, 236 Bulgaria blvd., Plovdiv 4003, BULGARIA}

\email{mmanev@pu.acad.bg, mmanev@yahoo.com}

%\thanks{The author was supported by the Matsumae International Foundation, Japan.}



\dedicatory{Dedicated to the 60th anniversary of Prof. Kouei SEKIGAWA}
\begin{abstract}
The tangent bundle as a $4n$-manifold is equipped with an almost hypercomplex
pseudo-Hermitian structure and it is characterized with respect to the relevant
classifications. A number of 8-dimensional examples of the considered type of
manifold are received from the known explicit examples in that manner.

\end{abstract}

%%% ----------------------------------------------------------------------
\maketitle
%%% ----------------------------------------------------------------------

%%%%%%%%%%%%%%%%%%%%%%%%%%%%%%%%%%%%%%%%%%%%%%%%%%%%%%%%%%%%%%%%%%%%%%%%%%%0
\section*{Introduction}

The geometry of the almost hypercomplex manifolds with Hermitian metric is
known (e.~g.~\cite{AlMa}). A parallel direction including indefinite metrics is
the developing of the geometry of the almost hypercomplex manifolds with
pseudo-Hermitian metric structure. It has a natural origination from the
geometry of the $n$-dimensional quaternionic Euclidean space.

The beginning was put by our joint works with K.~Gribachev and S.~Dimiev in
\cite{GrMaDi1} and \cite{GrMaDi2}. More precisely we have combined the ordinary
Hermitian metrics with the so-called by us \emph{skew-Hermitian metrics} with
respect to the almost complex structures of a hypercomplex structure.

First, we have ascertained four types of bilinear forms compatible with a
hypercomplex structure: three metrics and one K\"ahler form which are
skew-Hermitian or Hermitian with respect to the different almost complex
structures, i.~e. we have constructed a \emph{pseudo-Hermitian structure}.

Second, we have developed the notion of an \emph{almost hypercomplex
pseudo-Hermitian manifold}.

Third, we have considered mainly the (integrable) hypercomplex pseudo-Hermitian
manifolds and its special subclass the class of pseudo-hyper-K\"ahler
manifolds.

After that in \cite{Ma} we have continued the investigations as we have
constructed and characterized to a certain extent several 4-dimensional
examples of the considered manifolds.

The aim of the present work is to generate the considered structure on a
tangent bundle, to characterize it in terms of the corresponding manifolds and
to construct 8-dimensional examples relevant to the known 4-dimensional
examples.

In \secref{sec_ac} we recall the notions of the almost complex manifolds with
Hermitian metric or skew-Hermitian metrics and the almost hypercomplex
manifolds with pseudo-Hermitian metric structure.

In \secref{sec_TB} we remind the definitions and properties of the
prolongations of vector fields to tangent bundles known as their vertical and
horizontal lifts.

The main part of this paper is given in \secref{sec3}, where we consider an
almost skew-Hermitian manifold as a base manifold and we generate its tangent
bundle with a metric, which is a prolongation of the base metric, known as the
Sasaki metric. In that way we get the tangent bundle with almost hypercomplex
pseudo-Hermitian structure. We characterize differential-geometrically this
manifold having in mind also the relevant classifications of the received
manifolds.

The results of \secref{sec3} and the known examples from \cite{GrMaDi1} and
\cite{Ma} are used in the last \secref{sec4} to construct corresponding
8-dimensional examples in the discussed manner.

%%%%%%%%%%%%%%%%%%%%%%%%%%%%%%%%%%%%%%%%%%%%%%%%%%%%%%%%%%%%%%%%%%%%%%%%%%%%%%%%%1

\section{Some Kinds of Differentiable Manifolds with Almost Complex
Structures}\label{sec_ac}

\subsection{Almost Complex Manifolds with Hermitian Metric or Skew-Hermitian
Metric}

The notion of the \emph{almost complex manifold} $(M^{2n},J)$ is well-known.
There exists a possibility it to be equipped with two different kinds of
metrics. When $J$ acts as an isometry on each tangent space then the manifold
is an \emph{almost Hermitian manifold}. But, in the case when $J$ acts as an
anti-isometry on each tangent space, the notion of the so-called \emph{almost
skew-Hermitian manifold} is available. Let us consider more precisely the last
one.

Every $n$-dimensional complex Riemannian manifold induces a real
$2n$-dimensional manifold \\
$(M^{2n},J,g,\tilde g)$ with a complex structure $J$, a metric $g$ and an
associated metric $\tilde g = g (\cdot,J\cdot)$. Both metrics are indefinite of
signature $(n,n)$. This manifold we call an \emph{almost skew-Hermitian
manifold}. Such metric is known also as a B-metric or a Norden metric because
it is introduced by A.~P.~Norden in~\cite{No}. An almost skew-Hermitian
manifold is a \emph{skew-K\"ahler manifold} if $J$ is parallel with respect to
the Levi-Civita connection $\nabla$ of the metric $g$. The class of these
manifolds is contained in every other class of the almost skew-Hermitian
manifolds. A classification with respect to $\nabla J$ consisting of three
basic classes is given in~\cite{GaBo}. The class $\mathcal{W}_{1}$ is the main
class, i.~e. the class where $\nabla J$ has an explicit expression in terms of
the metric tensors and the Lie form.

%%%%%%%%%%%%%%%%%%%%%%%%%%%%%%%%%%%%%%%%%%%%%%%%%%%%%%%%%%%%%%%%%%%%%%%%%%%%%%%%%%%
\subsection{Almost Hypercomplex Manifolds with Pseudo-Hermitian Structure}~

Let us recall the notion of the almost hypercomplex structure $H$ on the
manifold $M^{4n}$. It is the triple $H=(J_\alpha)$ $(\alpha=1,2,3)$ of
anticommuting almost complex structures satisfying the property $J_3=J_1\circ
J_2$ (\cite{AlMa}, \cite{So}).

The introduction and the beginning of the study of the pseudo-Hermitian
structures on an almost hypercomplex manifold is given
in~\cite{GrMaDi1,GrMaDi2}. A pseudo-Riemannian metric $g$ of signature
$(2n,2n)$ on $(M^{4n},H)$ is introduced as follows
\begin{equation}\label{1}
g(\cdot,\cdot)=g(J_1\cdot,J_1\cdot)=-g(J_2\cdot,J_2\cdot)=-g(J_3\cdot,J_3\cdot).
\end{equation}
We called such metric a \emph{pseudo-Hermitian metric}. It generates a K\"ahler
2-form $\Phi$ and two pseudo-Hermitian metrics $g_2$ and $g_3$ by the following
way
\begin{equation}\label{G}
\Phi:=g(J_1\cdot,\cdot),\qquad g_2:=g(J_2\cdot,\cdot),\qquad
g_3:=g(J_3\cdot,\cdot).
\end{equation}
Let us note that $g$ ($g_2$, $g_3$, respectively) has an Hermitian
compatibility with respect to $J_1$ ($J_3$, $J_2$, respectively) and a
skew-Hermitian compatibility with respect to $J_2$ and $J_3$ ($J_1$ and $J_2$,
$J_1$ and $J_3$, respectively).

On the other hand, a quaternionic inner product $<\cdot,\cdot>$ in $\mathbb{H}$
generates in a natural way the bilinear forms $g$, $\Phi$, $g_2$ and $g_3$ by
the following decomposition: $<\cdot,\cdot>=-g+i\Phi+jg_2+kg_3$.

We called the structure $(H,G):=(J_1,J_2,J_3;g,\Phi,g_2,g_3)$ a
\emph{hypercomplex pseu\-do-Hermit\-ian structure} on $M^{4n}$ or shortly a
\emph{$(H,G)$-structure} on $M^{4n}$. We called the manifold $(M,H,G)$ an
\emph{almost hypercomplex pseudo-Hermitian manifold} or shortly an \emph{almost
$(H,G)$-manifold}.

It is well known, that the almost hypercomplex structure $H=(J_\alpha)$ is a
\emph{hypercomplex structure} if the Nijenhuis tensors
$%\[
N_\alpha(\cdot,\cdot)=
    \left[\cdot,\cdot \right]
    +J_\alpha\left[\cdot,J_\alpha \cdot \right]
    +J_\alpha\left[J_\alpha \cdot,\cdot \right]
    -\left[J_\alpha \cdot,J_\alpha \cdot \right]
$ %\]
vanish for each $\alpha=1,2,3$. Moreover, one $H$ is hypercomplex iff two
of $N_\alpha$ vanish.

We  introduced \emph{structural $(0,3)$-tensors} of the almost $(H,G)$-manifold
by
\[
F_\alpha (x,y,z)=g\bigl( \left( \nabla_x J_\alpha
\right)y,z\bigr)=\bigl(\nabla_x g_\alpha\bigr) \left( y,z \right) \quad
(\alpha=1,2,3),
\]
where $\nabla$ is the Levi-Civita connection generated by $g$ and $x, y, z \in
T_pM$ at any $p\in M$.

There are valid relations between $F_\alpha$'s, e.~g.
$F_1(\cdot,\cdot,\cdot)=F_2(\cdot,J_3\cdot,\cdot)+F_3(\cdot,\cdot,J_2\cdot)$.

Since $g$ is a Hermitian metric with respect to $J_1$, according to
Gray-Hervella \cite{GrHe} the basic subclass $\W_4$ of the Hermitian manifolds
is determined by
\[
\begin{array}{c}
F_1(x,y,z)=\frac{1}{2(2n-1)} \left[ g(x,y)\theta_1(z)-g(x,z)\theta_1(y)
-g(x,J_1y)\theta_1(J_1z)+g(x,J_1z)\theta_1(J_1y) \right],
\end{array}
\]
where $\theta_1(\cdot)=g^{ij}F_1(e_i,e_j,\cdot)=\delta\Phi(\cdot)$ is the
\emph{Lie form} for any basis $\{e_i\}_{i=1}^{4n}$, and $\delta$ -- the
coderivative.

On other side, the metric $g$ is a skew-Hermitian one with respect to $J_2$ and
$J_3$. According to the classification of the almost complex manifolds with
skew-Hermitian metric (Norden metric or B-metric) given in \cite{GaBo} the
basic classes are defined as follows:
\begin{equation}\label{GaBo}
\begin{array}{c}
\mathcal{W}_1: F_\alpha(x,y,z)=\frac{1}{4n} \left[
g(x,y)\theta_\alpha(z)+g(x,z)\theta_\alpha(y) +g(x,J_\alpha
y)\theta_\alpha(J_\alpha z)
    +g(x,J_\alpha z)\theta_\alpha(J_\alpha y)\right],\\[6pt]
\mathcal{W}_2: \mathop{\makebox{\huge $\sigma$}} \limits_{x,y,z}
F_\alpha(x,y,J_\alpha z)=0,\qquad
\mathcal{W}_3: \mathop{\makebox{\huge $\sigma$}} \limits_{x,y,z}
F_\alpha(x,y,z)=0,
\end{array}
\end{equation}
where \(\theta_\alpha(\cdot)=g^{ij}F_\alpha (e_i,e_j,\cdot)\), \(\alpha =2,3\),
are the corresponding Lie forms for an arbitrary basis \(\{e_i\}_{i=1}^{4n}\).

Here, we denote the main subclasses of the respective complex manifolds by
$\W(J_\alpha)$, where $\W(J_1):=\W_4(J_1)$ in~\cite{GrHe}, and
$\W(J_\alpha):=\W_1(J_\alpha)$ for $\alpha=2,3$ in~\cite{GaBo}.

We obtained a sufficient condition an almost $(H,G)$-manifold to be an
integrable one:

\begin{thm}[\cite{GrMaDi1}] \label{t31}
Let $(M,H,G)$ belongs to $ \W(J_\alpha) \bigcap \W(J_\beta)$. Then $(M,H,G)$ is
of class $\W(J_\gamma)$ for all cyclic permutations $(\alpha, \beta, \gamma)$
of $(1,2,3)$.
\end{thm}

We say that a pseudo-Hermitian manifold is a \emph{pseudo-hyper-K\"ahler
manifold} and denote $(M,H,G)$ $\in \K$, if $F_\alpha=0$ for every
$\alpha=1,2,3$, i.~e. the manifold is of K\"ahlerian type with respect to each
$J_\alpha$.

\begin{thm}[\cite{GrMaDi1}]\label{t33}
If $(M,H,G) \in \K(J_\alpha)\bigcap \W(J_\beta)$ $(\alpha\neq\beta \in
\{1,2,3\})$ then $(M,H,G) \in \K$ .
\end{thm}

We gave a geometric characteristic of the pseudo-hyper-K\"ahler manifolds
according to the curvature tensor $R=[\nabla ,\nabla ] - \nabla_{[\ ,\ ]}$
induced by the Levi-Civita connection.
\begin{thm}[\cite{GrMaDi2}]\label{R=0}
Each pseudo-hyper-K\"ahler manifold is a flat pseudo-Riemann\-ian manifold with
signature $(2n,2n)$.
\end{thm}

%%%%%%%%%%%%%%%%%%%%%%%%%%%%%%%%%%%%%%%%%%%%%%%%%%%
\section{The Tangent Bundle}\label{sec_TB}

Let us consider a $2n$-dimensional differentiable manifold $M$ and then $T_p M$
is the tangent space at a point $p$ of $M$. It is known that the set
\[
TM=\bigcup_{p\in M}T_p M
\]
is called the \emph{tangent bundle} over the manifold $M$~(\cite{St}).

For any point $u$ of $TM$ such that $u\in T_p M$ the correspondence
$u\rightarrow p$ defines the bundle projection $\pi: TM \rightarrow M$ by
$\pi(u)=p$. The set $\pi^{-1}(p)$, i.~e. $T_p M$, is called the \emph{fiber} of
$TM$ over $p\in M$ and $M$ -- the \emph{base manifold} of $TM$. The base space
$M$ is a submanifold differentiably imbedded in $TM$.

Let $(x^i)$ are the local coordinates of a point $p$ in a neighborhood $U
\subset M$. Then a point $u\in T_p M$ with coordinates $(y^i)$ in $T_p M$ with
respect to the natural basis $\left\{\frac{\partial}{\partial x^i}\right\}$ is
represented by the ordered pair $(p,u)$. Therefore an arbitrary point
$\tilde{p}$ in $TM$ has local coordinates of the kind $(x^i,y^i)$ in the open
set $\pi^{-1}(U)\subset TM$. These coordinates is called the \emph{induced
coordinates} in $\pi^{-1}(U)$ from $(x^i)$.

Obviously, the tangent bundle of a base manifold $M$ with almost complex
structure is of dimension $4n$.

Let us denote by $\mathcal{T}^r_s(M)$ the set of all differentiable tensor
fields of type $(r,s)$ in $M$ and let
$\mathcal{T}(M)=\sum_{r,s=0}^\infty\mathcal{T}^r_s(M)$ be the set of all tensor
fields in $M$. Similarly, let us denote by $\mathcal{T}^r_s(TM)$ and
$\mathcal{T}(TM)$ respectively the corresponding sets of tensor fields in the
tangent bundle $TM$.

\subsection{Vertical Lifts}
The \emph{vertical lift $f^V$ of a function $f$} in $M$ is called the
composition of $\pi: TM \rightarrow M$ and $f: M \rightarrow \mathbb{R}$, i.~e.
%\begin{equation}
$ f^V=f\circ \pi, $
%\end{equation}
where $f \in \mathcal{T}^0_0(M)$ and $f^V\in \mathcal{T}^0_0(TM)$.

Since a correspondence $TM \rightarrow \mathbb{R}$ is possible to be considered
as the action of a 1-form in $M$ or a function in $TM$, then if $\omega$ is a
1-form in $M$, it is regarded, in a natural way, as a function in $TM$, denoted
by $\imath\omega$.

Let $\tilde{X}\in\mathcal{T}^1_0(TM)$ be such that $\tilde{X}f^V=0$ for all $f
\in \mathcal{T}^0_0(M)$. Then $\tilde{X}$ is called a \emph{vertical vector
field} in $TM$. The vector field $X^V$ in $TM$ defined by
$X^V(\imath\omega)=\left(\omega(X)\right)^V$ for an arbitrary 1-form $\omega$
in $M$, is called the \emph{vertical lift of vector field} $X$ in $M$ to $TM$.
If $X$ has components $X^i$ in $M$, then $X^V$ has components
\begin{equation}\label{XV-ind}
X^V:
\left(%
\begin{array}{c}
  0 \\
  X^i \\
\end{array}%
\right)
\end{equation}
 with respect to the induced coordinates in $TM$.
Consequently, it is clear, that \( \left(\frac{\partial}{\partial
x^i}\right)^V=\frac{\partial}{\partial y^i}\).

There are valid the following formulas for any $X, Y\in\mathcal{T}^1_0(M)$ and
$f\in\mathcal{T}^0_0(M)$.
\[%\begin{equation}
X^Vf^V=0,\quad
\left(X+Y\right)^V=X^V+Y^V,\quad \left(fX\right)^V=f^VX^V,\quad
\left[X^V,Y^V\right]=0.
\]%\end{equation}

%%%%%%%%%%%%%%%%%%%%---------------------------------
\subsection{Horizontal Lifts}

Let $\nabla$ be the Levi-Civita connection in $M$ induced by the
pseudo-Riemann\-ian metric $g$. \emph{Horizontal lift $f^H$ of a function} $f$
in
$M$ to $TM$ is defined by $%\[
f^H=\imath(df)-\gamma(\nabla f),$ %\]
where $\gamma(\nabla f)=y^k\nabla_k f$.
There follows immediately that $f^H=0$ for any $f \in \mathcal{T}^0_0(M)$.

Suppose that $X$ and $\nabla$ have local components $X^k$ and $\Gamma_{ij}^k$,
respectively, in $M$. \emph{Horizontal lift $X^H$ of a vector field} $X$ in $M$
to $TM$ is called the vector field in $TM$ with components
\begin{equation}\label{XH-ind}
X^H: \left(%
\begin{array}{c}
  X^k \\
  -y^i\Gamma_{ij}^kX^j \\
\end{array}%
\right)
\end{equation}
with respect to the induced coordinates in $TM$.

For the action of $X^H$ we have that
\(
    X^Hf^V=(Xf)^V,\ %\quad
        \omega^V(X^H)=\{\omega(X)\}^V,\ %\quad
            J^VX^H=(JX)^V.
\)

%%%%%%%%%%%%%%%%%%%%%%%%%%%%%%%%%%%%%%%%%%%%%%%%%%%%%%%%%%%%
\subsection{Adapted Frames}
In each coordinate neighborhood $\{U,x^k\}$ of $M$, $\dim M=2n$, we denote
$e_{(i)}=\frac{\partial}{\partial x^i}$. Then $4n$ local vector fields
$e_{(i)}^H$ and $e_{(i)}^V$ form a basis of the tangent space
$T_p\left(TM\right)$ at each point $\tilde{p}\in \pi^{-1}(p)$ and their
components are given by
\[
e_{(i)}^H:
\left(%
\begin{array}{c}
  \delta_i^k \\
  -y^j\Gamma_{ij}^k \\
\end{array}%
\right) , \qquad e_{(i)}^V:
\left(%
\begin{array}{c}
  0 \\
  \delta_i^k \\
\end{array}%
\right)
\]
with respect to the induced coordinates $(x^k,y^k)$ in $TM$. The set
$\left\{e_{(i)}^H,e_{(i)}^V\right\}$ is called the \emph{frame adapted to the
connection} $\nabla$ in $\pi^{-1}(U)$. On putting $\tilde{e}_{(i)}=e_{(i)}^H,
\tilde{e}_{(\bar{i})}=e_{(i)}^V$ we write the adapted frame as
$\{\tilde{e}_{(A)}\}=\{\tilde{e}_{(i)},\tilde{e}_{(\bar{i})}\}$.

Then the horizontal and vertical lifts of $X$ %($\omega$, resp.)
with local components $X^k$ %and $\omega_k$ have,
are components
\[%\begin{equation}\label{XHV}
X^H:
\left(%
\begin{array}{c}
  X^k \\
  0 \\
\end{array}%
\right) , \qquad X^V:
\left(%
\begin{array}{c}
  0 \\
  X^k \\
\end{array}%
\right)
\]%\end{equation}
with respect to the adapted frame $\{\tilde{e}_{(A)}\}$ in $TM$.

%%%%%%%%%%%%%%%%%%%%%%%%%%%%%%%%%%%%%%%%%%%%%%%%%%%%%%%%%%%%%%%%%%%%%%%%%%%%%%%%%
\section{The Tangent Bundle with an Almost Hypercomplex Pseudo-Hermitian
Structure}\label{sec3}

Our purpose is a determination of an almost hypercomplex pseudo-Hermitian
structure $(H,G)$ on $TM$ when the base manifold $M$ has an almost
skew-Hermitian structure $(J,g,\tilde{g})$.

We will use the horizontal and vertical lifts of the vector fields on $M$ to
get the corresponding components of the considered tensor fields on $TM$. These
components are sufficient to describe the characteristic tensor fields on $TM$
in general.

%%%%%%%%%%%%%%%%%%%%%%%%%%%%%%%%%%%%%%%%%%%%%%%%%%%%%%%%%%%%%%%%%%%%%%%%%%%%%%
\subsection{An Almost Hypercomplex Structure on the Tangent Bundle}

As it is known \cite{YaIs}, for any affine connection in $M$, the induced
horizontal and vertical distributions in $TM$ are mutually complementary. Then
we define tensor field $J_1$, $J_2$ and $J_3$ in $TM$ by their action over the
horizontal and vertical lifts of an arbitrary vector field in $M$ as follows:
\begin{equation}\label{H}
J_1:\ X^H\rightarrow X^V,\; X^V\rightarrow -X^H;\quad J_2:\ X^H\rightarrow
(JX)^V,\; X^V\rightarrow (JX)^H;\quad J_3:=J_1\circ J_2, %\ X^H\rightarrow -(JX)^H,\;
%X^V\rightarrow (JX)^V,
\end{equation}
where $J$ is the given almost complex structure in $M$.

By direct computations we get the following

\begin{prop}
There exists an almost hypercomplex structure $H=(J_1,J_2,J_3)$, defined
by~\eqref{H} in $TM$ over an almost complex manifold $(M,J)$ with an affine
connection $\nabla$. The received $4n$-dimensional manifold is an almost
hypercomplex manifold $(TM,H)$.
\end{prop}

Let $N_\alpha$ denotes the Nijenhuis tensor regarding $J_\alpha$ for each
$\alpha=1,2,3$ and $\tilde{X},\tilde{Y}\in \mathcal{T}^0_1(TM)$, i.~e.
\begin{equation}\label{Na-def}
N_\alpha(\tilde{X},\tilde{Y})=[\tilde{X},\tilde{Y}]
+J_\alpha[J_\alpha\tilde{X},\tilde{Y}] +J_\alpha[\tilde{X},J_\alpha\tilde{Y}]
-[J_\alpha\tilde{X},J_\alpha\tilde{Y}].
\end{equation}

Having in mind \eqref{XV-ind} and \eqref{XH-ind} we receive by direct
computations the following
\begin{lem}\label{l[]}
Let $\nabla$ be a torsion-free affine connection in $M$ and $R$ be its
curvature tensor. Then for any $X, Y\in \mathcal{T}^0_1(M)$ at $u\in T_p M$ we
have
\[
    \begin{array}{ll}
    [X^H,Y^H]=[X,Y]^H-\{R(X,Y)u\}^V,\quad
    &[X^H,Y^V]=\left(\nabla_XY\right)^V,\\[4pt]
    [X^V,Y^V]=0,
        &[X^V,Y^H]=-\left(\nabla_YX\right)^V.
   \end{array}
\]
\end{lem}

Using the definitions \eqref{H}, \eqref{Na-def} and \lemref{l[]} we get
\begin{prop}\label{Na}
Let $(M,J)$ be an almost complex manifold. Then the Nijenhuis tensors of the
structure $H$ in $TM$ for the horizontal and vertical lifts have the form
\begin{equation}\label{N_1}
\begin{array}{l}
-N_1(X^H,Y^H)=N_1(X^V,Y^V)=\left\{R(X,Y)u\right\}^V,\;\\[4pt]
N_1(X^H,Y^V)=N_1(X^V,Y^H)=-\left\{R(X,Y)u\right\}^H;\\[4pt]
\end{array}
\end{equation}
\begin{equation}\label{N_2}
\begin{array}{l}
N_2(X^H,Y^H)=\left\{J(\nabla_X J)(Y)-J(\nabla_Y
J)(X)\right\}^H-\left\{R(X,Y)u\right\}^V,\;\\[4pt]
N_2(X^V,Y^V)=-\left\{(\nabla_{JX}J)(Y)+(\nabla_{JY}
J)(X)\right\}^H+\left\{R(JX,JY)u\right\}^V,\; \\[4pt]
N_2(X^H,Y^V)=\left\{J(\nabla_X
J)(Y)+(\nabla_{JY} J)(X)\right\}^V-\left\{JR(X,JY)u\right\}^H,\;\\[4pt]
N_2(X^V,Y^H)=-\left\{(\nabla_{JX}
J)(Y)+(J\nabla_{Y} J)(X)\right\}^V-\left\{JR(JX,Y)u\right\}^H;\\[4pt]
\end{array}
\end{equation}
\begin{equation}\label{N_3}
\begin{array}{l}
N_3(X^H,Y^H)=\left\{-R(X,Y)u+R(JX,JY)u+JR(JX,Y)u+JR(X,JY)u\right\}^V, \\[4pt]
%\phantom{N_3(X^H,Y^H)=}+\left\{N(X,Y)\right\}^H,\qquad\qquad
%
N_3(X^H,Y^V)=\left\{(\nabla_{JX}
J)(Y)-(\nabla_{X} J)(JY)\right\}^V,\;\\[4pt]
N_3(X^V,Y^H)=-\left\{(\nabla_{Y} J)(JX)-(\nabla_{JY} J)(X)\right\}^V,\quad
N_3(X^V,Y^V)=0.
\end{array}
\end{equation}
for any $X, Y\in \mathcal{T}^0_1(M)$ at $u\in T_p M$.
\end{prop}

The last equalities for $N_\alpha$ imply the following necessary and sufficient
conditions for integrability of $J_\alpha$ and $H$.
\begin{thm}\label{tH}\  Let $TM$ be the tangent bundle manifold
with an almost hypercomplex structure $H=(J_1,J_2,J_3)$ defined as in~\eqref{H}
and $M$ be its base manifold with an almost complex structure $J$. Then the
following inter\-connections hold:
\begin{enumerate}
    \item $(TM,J_1)$ is
complex iff $M$ is flat; \emph{(see also \cite{To})}
    \item $(TM,J_\beta)$ for $\beta=2$
or $3$ is complex iff $M$ is flat and $J$ is parallel;
    \item $(TM,H)$ is hypercomplex iff $M$ is flat and $J$ is parallel.
\end{enumerate}
\end{thm}
\begin{cor}\
\begin{enumerate}
    \item $(TM,J_2)$ is complex iff $(TM,J_3)$ is complex.
    \item If $(TM,J_2)$ or $(TM,J_3)$ is complex then $(TM,H)$ is hypercomplex.
\end{enumerate}
\end{cor}

%%%%%%%%%%%%%%%%%%%%%%%%%%%%%%%%%%%%%%%%%%%%%%%%%%%%%%%%%%%%%%%%%%%%%%%%%%%%%%%%%%%%%%
\subsection{The Sasaki Metric on the Tangent Bundle}

Let us introduce the \emph{Sasaki metric} $\hat{g}$ on $TM$ as in~\cite{Sa}
(i.~e. $\hat{g}$ is the so-called \emph{diagonal lift} of the base metric $g$)
defined by
\begin{equation}\label{5g}
\hat{g}(X^H,Y^H)=\hat{g}(X^V,Y^V)=g(X,Y),\quad \hat{g}(X^H,Y^V)=0.
\end{equation}

Having in mind the definition~\eqref{H} of the structure $H$, we verify
immediately that the Sasaki metric satisfies the properties~\eqref{1} and
therefore it is valid the following

\begin{thm}\label{Sas}
The tangent bundle $TM$ equipped with the almost hypercomplex structure $H$ and
the Sasaki metric $\hat{g}$, defined by~\eqref{H} and~\eqref{5g}, respectively,
is an almost hypercomplex pseudo-Hermitian manifold denoting by
$(TM,H,\hat{G})$.
\end{thm}

Since $\hat{\nabla}$ is the Levi-Civita connection of $\hat{g}$ on $TM$ as
$\nabla$ of $g$ on $M$, then using the property of type
\[
2g(\nabla_X Y,Z)=Xg(Y,Z)+Y g(X,Z)-Z g(X,Y)+g([X,Y],Z)+g([Z,X],Y)+g([Z,Y],X)
\]
we obtain the covariant derivatives of the horizontal and vertical lifts of
vector fields on $TM$ as follows

\begin{lem}\label{nabli} For any $X, Y\in \mathcal{T}^0_1(M)$ at $u\in T_p M$
\[
\begin{array}{ll}
\hat{\nabla}_{X^H}Y^H=(\nabla_XY)^H-\frac{1}{2}\{R(X,Y)u\}^V,\quad &
\hat{\nabla}_{X^H}Y^V=\frac{1}{2}\{R(u,Y)X\}^H+(\nabla_XY)^V,\\[4pt]
\hat{\nabla}_{X^V}Y^H=\frac{1}{2}\{R(u,X)Y\}^H,
& \hat{\nabla}_{X^V}Y^V=0.
\end{array}
\]
\end{lem}

After that we calculate the components of the curvature tensor $\hat{R}$ of
$\hat{\nabla}$ with respect to the horizontal and vertical lifts of the vector
fields on $M$ and we receive
\begin{prop}\label{R}
The following interconnections between the curvature tensors $R$ and $\hat{R}$
corresponding to the metrics $g$ and $\hat{g}$ are valid
\[
\begin{array}{l}
\hat{R}(X^H,Y^H,Z^H,W^H)=R(X,Y,Z,W)
+\frac{1}{4}\{g\left(R(W,X)u,R(Y,Z)u\right)-g\left(R(W,Y)u,R(X,Z)u\right)\}\\[4pt]
\phantom{\hat{R}(X^H,Y^H,Z^H,W^H)=}
-\frac{1}{2}g\left(R(X,Y)u,R(Z,W)u\right),\\[4pt]
\hat{R}(X^H,Y^H,Z^H,W^V)=-\frac{1}{2}g\left(
(\nabla_XR)(Y,Z)u-(\nabla_YR)(X,Z)u,W
\right),\\[4pt]
\hat{R}(X^H,Y^H,Z^V,W^V)=R(X,Y,Z,W)-\frac{1}{4}\{g\left(R(u,W)X,R(u,Z)Y\right)
-g\left(R(u,W)Y,R(u,Z)X\right)\},\\[4pt]
\hat{R}(X^H,Y^V,Z^H,W^H)=\frac{1}{2}g\left( (\nabla_XR)(u,Y)Z,W \right),\\[4pt]
\hat{R}(X^H,Y^V,Z^H,W^V)=\frac{1}{2}R(X,Z,Y,W)-\frac{1}{4}g\left(R(u,Y)Z,R(u,W)X\right),\\[4pt]
\hat{R}(X^V,Y^V,Z^H,W^H)=R(X,Y,Z,W)
-\frac{1}{4}\{g\left(R(u,Y)Z,R(u,X)W\right)-g\left(R(u,X)Z,R(u,Y)W\right)\},\\[4pt]
\hat{R}(X^V,Y^V,Z^H,W^V)=\hat{R}(X^H,Y^V,Z^V,W^V)=\hat{R}(X^V,Y^V,Z^V,W^V)=0.
\end{array}
\]
for the lifts $(\cdot)^H, (\cdot)^V \in \mathcal{T}^0_1(TM)$ of any $X, Y, Z, W
\in \mathcal{T}^0_1(M)$ at $u\in T_p M$.
\end{prop}

Hence we get
\begin{thm}\label{flat}
The tangent bundle $(TM,\hat{g})$ is flat if and only if its base
manifold $(M,g)$ is flat.
\end{thm}

\begin{rem}
The results of the last three statements \ref{nabli}--\ref{flat}
are confirmed also by \cite{Ko} (see also \cite{GuKa}), where $g$
is a Riemannian metric.
\end{rem}

According to \cite{Ko}, the tangent bundle with Sasaki metric of a
Riemannian manifold is never locally symmetric unless the base
manifold is locally Euclidean. Having in mind the last reasonings,
the proof of the corresponding theorem in \cite{Ko} and that $g$
is a pseudo-Riemannian metric in our case, we obtain the following

\begin{concl}
If the tangent bundle $(TM,\hat{g})$ is locally symmetric then the
curvature tensor $R$ of its base pseudo-Riemannian manifold
$(M,g)$ is zero or isotropic, i.e. $R=0$ or $g(R,R)=0, R\neq 0$).
\end{concl}

%%%%%%%%%%%%%%%%%%%%%%%%%%%%%%%%%%%%%%%%%%%%%%%%%%%%%%%%%%%%%%%%%%%%%%%%%%%%%%%%%
\subsection{The Tangent Bundles with Almost Hermitian, Almost Skew-Hermitian
and Almost Hypercomplex Pseudo-Hermitian Structures}

Suppose that $(M,J)$ is an almost complex manifold with skew-Hermitian metrics
$g$, $\tilde{g}$ and that $(TM,H)$ is its almost hypercomplex tangent bundle
with the pseudo-Hermitian metric structure
$\hat{G}=(\hat{g},\hat{\Phi},\hat{g}_2,\hat{g}_3)$ derived (as in~\eqref{G})
from the Sasaki metric $\hat{g}$ on $TM$ - the diagonal lift of $g$. The
generated $4n$-dimensional manifold we will denote by $(TM,H,\hat{G})$.

To characterize the structural tensors
$F_\alpha(\cdot,\cdot,\cdot)=\hat{g}\left((\hat{\nabla}_{\cdot}
J_\alpha)(\cdot),\cdot\right)$ at each $u\in T_pM$ on $(TM,H,\hat{G})$ we use
\lemref{nabli} and the definitions \eqref{H} of $J_\alpha$, whence we obtain
the following

\begin{prop}\label{Fa} The nonzero components of $F_\alpha$ with respect to
the lifts of the vector fields depend on structural tensor $F$ and the
curvature tensor $R$ on $(M,J,g,\tilde{g})$ by the following way:
\begin{equation}\label{5F1}
\begin{array}{l}
-F_1(X^H,Y^H,Z^H)=F_1(X^H,Y^V,Z^V)=F_1(X^V,Y^H,Z^V)=F_1(X^V,Y^V,Z^H)\\[4pt]
\phantom{-F_1(X^H,Y^H,Z^H)}=\frac{1}{2}R(Y,Z,X,u);
\end{array}
\end{equation}
\begin{equation}\label{5F2}
\begin{array}{l}
F_2(X^H,Y^H,Z^H)=-\frac{1}{2}R(X,Y,JZ,u)+\frac{1}{2}R(Z,X,JY,u),\\[4pt]
F_2(X^H,Y^V,Z^V)=\frac{1}{2}R(X,JY,Z,u)-\frac{1}{2}R(JZ,X,Y,u),\\[4pt]
F_2(X^H,Y^H,Z^V)=F_2(X^H,Y^V,Z^H)=F(X,Y,Z),\\[4pt]
F_2(X^V,Y^H,Z^V)=\frac{1}{2}R(Y,JZ,X,u),\quad
F_2(X^V,Y^V,Z^H)=-\frac{1}{2}R(JY,Z,X,u);
\end{array}
\end{equation}
\begin{equation}\label{5F3}
\begin{array}{l}
F_3(X^H,Y^H,Z^H)=-F_3(X^H,Y^V,Z^V)=-F(X,Y,Z),\quad\\[4pt]
F_3(X^H,Y^H,Z^V)=-\frac{1}{2}R(X,JY,Z,u)-\frac{1}{2}R(X,Y,JZ,u),\\[4pt]
F_3(X^H,Y^V,Z^H)=\frac{1}{2}R(Z,X,JY,u)+\frac{1}{2}R(JZ,X,Y,u),\\[4pt]
F_3(X^V,Y^H,Z^H)=\frac{1}{2}R(JY,Z,X,u)-\frac{1}{2}R(Y,JZ,X,u).\\[4pt]
\end{array}
\end{equation}
where $(\cdot)^H, (\cdot)^V \in \mathcal{T}^0_1(TM)$ are the lifts of any $X,
Y, Z \in \mathcal{T}^0_1(M)$ at $u\in T_p M$.
\end{prop}

Hence we compute the corresponding Lie forms with respect to the adapted frame.

Let $\{e_i^H,e^{H}_{\bar{i}},e_i^V,e^{V}_{\bar{i}}\}$ be the adapted frame of
type $(+,\dots,+,-,\dots,-,+,\dots,+,-,\dots,-)$ at each point of $TM$ derived
by the orthonormal basis $\{e_{(i)},e_{(\bar{i})}\}$ of signature $(n,n)$ at
each point of $M$ denoting $e_{(\bar{i})}=Je_i$. The indices $i$ and $\bar{i}$
run over the ranges $\{1,\dots,n\}$ and $\{n+1,\dots,2n\}$, respectively. For
example
\[
\begin{array}{l}
\theta_3(Z^H)=\sum_{i=1}^n\left\{F_3(e_i^H,e_i^H,Z^H)-F_3(e_{\bar{i}}^H,e_{\bar{i}}^H,Z^H)
+F_3(e_i^V,e_i^V,Z^H)-F_3(e_{\bar{i}}^V,e_{\bar{i}}^V,Z^H)\right\}\\[4pt]
\phantom{\theta_3(Z^H)}
=\sum_{i=1}^n\left\{-F(e_i,e_i,Z)+F(e_{\bar{i}},e_{\bar{i}},Z)\right\}=-\theta(Z),\\[4pt]
\theta_3(Z^V)=\sum_{i=1}^n\left\{F_3(e_i^H,e_i^H,Z^V)-F_3(e_{\bar{i}}^H,e_{\bar{i}}^H,Z^V)
+F_3(e_i^V,e_i^V,Z^V)-F_3(e_{\bar{i}}^V,e_{\bar{i}}^V,Z^V)\right\}\\[4pt]
\phantom{\theta_3(Z^V)}
=\sum_{i=1}^n\left\{-\frac{1}{2}R(e_i,e_{\bar{i}},Z,u)
-\frac{1}{2}R(e_{\bar{i}},e_i,Z,u)\right\}=0.\\[4pt]
\end{array}
\]
Then we have
\begin{prop}~
\begin{enumerate}
    \item $(TM,J_1,\hat{g})$ has a zero Lie form $\theta_1$;
    \item $(TM,J_2,\hat{g})$ has a zero Lie form $\theta_2$ iff
        $(M,J,g,\tilde{g})$ has a zero Lie form $\theta$ and a zero
        associated Ricci tensor $\tilde{\rho}$;
    \item $(TM,J_3,\hat{g})$ has a zero Lie form $\theta_3$ iff
        $(M,J,g,\tilde{g})$ has a zero Lie form $\theta$.
\end{enumerate}
\end{prop}

According (\ref{5F1}) the following property is valid.
\begin{thm}
$(TM,J_1,\hat{g})$ is an almost K\"ahler manifold and it is a K\"ahler manifold
if and only if $(M,J,g,\tilde{g})$ is flat.
\end{thm}
\begin{rem}
By comparison with the Riemannian case, in \cite{To} it is shown
that $(TM,J_1,\hat{g})$ is almost K\"ahlerian (i.e. symplectic)
for any Riemannian metric $g$ on the base manifold when the
connection used to define the horizontal lifts is the Levi-Civita
connection.
\end{rem}

Having in mind (\ref{5F1})--(\ref{5F3}) and \thmref{tH} it is easy
to conclude the following sequel of properties:
\begin{prop}\
\begin{enumerate}
    \item $(TM,J_1,\hat{g})$ is K\"ahlerian iff $M$ is flat.
    \item $(TM,J_\beta,\hat{g})$ for $\beta=2$ or $3$ is
        skew-K\"ahlerian iff $(M,J,g,\tilde{g})$ is flat and
        skew-K\"ahlerian.
    \item $(TM,H,G)$ is a pseudo-hyper-K\"ahler manifold iff
        $(M,J,g,\tilde{g})$ is flat and skew-K\"ahlerian.
\end{enumerate}
\end{prop}
\begin{cor}\
\begin{enumerate}
    \item $(TM,J_2,\hat{g})$ is skew-K\"ahlerian iff
        $(TM,J_3,\hat{g})$ is skew-K\"ahlerian.
    \item If
        $(TM,J_2,\hat{g})$ or $(TM,J_3,\hat{g})$ is skew-K\"ahlerian then
        $(TM,H,\hat{G})$ is pseudo-hyper-K\"ahlerian.
\end{enumerate}
\end{cor}
\begin{cor}\
\begin{enumerate}
    \item The only complex manifolds $(TM,J_\alpha,\hat{g})$ for
        some $\alpha=1,2,3$ are the K\"ahler
        manifolds with respect to that $J_\alpha$.
    \item The only hypercomplex manifolds $(TM,H,\hat{G})$ are the pseudo-hyper-K\"ahler
        manifolds.
\end{enumerate}
\end{cor}

Let us recall that the class $\left\{\W_2 \oplus \W_3\right\}$ of the almost
skew-Hermitian manifolds is determined by the condition for vanishing of the
Lie form, the class $\W_3$ by $\mathop{\makebox{\huge$\sigma$}} F=0$ and the
class $\W_0$ by $F=0$ (i.~e. the last is the class of the skew-K\"ahler
manifolds). Then we obtain the following
\begin{thm}Let $M$ be $(M,J,g,\tilde{g})$ and $TM$ be $(TM,H,\hat{G})$.
\begin{enumerate}
    \item $TM\in\left\{\W_2 \oplus \W_3\right\}(J_2)$ iff
        $M\in\left\{\W_2 \oplus \W_3\right\}(J)$ and
        $\rho=\tilde{\rho}=0$;
    \item $TM\in\W_3(J_2)$ iff $M\in\W_0(J)$ and
        $\rho=\tilde{\rho}=0$;
    \item $TM\in\left\{\W_2 \oplus \W_3\right\}(J_3)$ iff
        $M\in\left\{\W_2 \oplus \W_3\right\}(J)$;
    \item $TM\in\W_3(J_3)$ iff
        $M\in\W_0(J)$,
\end{enumerate}
where $\rho$ and $\tilde{\rho}$ denote the Ricci tensor and the associated
Ricci tensor on $M$, respectively.
\end{thm}

If we set some additional conditions for $(M,J,g,\tilde{g})$ then we receive
the corresponding specialization of $(TM,H,\hat{G})$ given in the next
properties:
\begin{prop} Let $M$ be $(M,J,g,\tilde{g})$ and $TM$ be $(TM,H,\hat{G})$.
\begin{enumerate}
    \item If $M\in \left\{\W_2 \oplus \W_3\right\}(J)$ then
        $TM\in \AK(J_1)\cap \left\{\W_1\oplus\W_2\oplus
        \W_3\right\}(J_2)\cap \left\{\W_2\oplus \W_3\right\}(J_3)$.
    \item If $M\in \left\{\W_2 \oplus \W_3\right\}(J)$ and
        $\rho=\tilde{\rho}=0$ then\\
        \phantom{If $M\in \left\{\W_2 \oplus \W_3\right\}(J)$ and
        $\rho=\tilde{\rho}=0$}
        $TM\in \AK(J_1)\cap
        \left\{\W_2\oplus \W_3\right\}(J_2)\cap \left\{\W_2\oplus
        \W_3\right\}(J_3)$.
    \item If $M\in \left\{\W_2 \oplus \W_3\right\}(J)$ and $R=0$
        then $TM\in \K(J_1)\cap \left\{\W_2\oplus
        \W_3\right\}(J_2)\cap \left\{\W_2\oplus \W_3\right\}(J_3)$.
\end{enumerate}
\end{prop}
\begin{prop} Let $M$ be $(M,J,g,\tilde{g})$ and $TM$ be $(TM,H,\hat{G})$.
\begin{enumerate}
    \item If $M\in \W_0(J)$ then $TM\in \AK(J_1)\cap
        \left\{\W_1\oplus\W_2\oplus \W_3\right\}(J_2)\cap \W_3(J_3)$.
    \item If $M\in \W_0(J)$ and $\rho=\tilde{\rho}=0$ then
        $TM\in \AK(J_1)\cap \W_3(J_2)\cap \W_3(J_3)$.
    \item If $M\in \W_0(J)$ and $R=0$ then $TM\in
        \K(J_1)\cap \W_0(J_2)\cap \W_0(J_3)$.
\end{enumerate}
\end{prop}

In above properties $\AK(J_1)$ and $\K(J_1)$ denote the classes regarding $J_1$
of the almost K\"ahler manifolds and the K\"ahler manifolds determined by the
conditions $\theta_1=0$ and $F_1=0$, respectively.
%%%%%%%%%%%%%%%%%%%%%%%%%%%%%%%%%%%%%%%%%%%%%%%%%%%%%%%%%%%%%%%%%%%%%%%%%%%%%%%%%%2

\section{Eight-Dimensional Derived Examples}\label{sec4}

Let us recall the developed examples of 4-dimensional almost hypercomplex
pseudo-Hermitian manifolds. In \cite{GrMaDi1} and \cite{Ma} we have constructed
a lot of explicit examples of the investigated manifolds.

At first in \cite{GrMaDi1}, it is considered a pseudo-Riemannian spherical
manifold $S_2^4$ in pseudo-Euclid\-ean vector space $\mathbb{R}_2^5$ of
signature $(2,3)$. It admits a hypercomplex pseudo-Hermitian structure, with
respect to which it is of the class $\W=\bigcap_{\alpha=1}^{3}\W(J_\alpha)$.

Secondly, in the same paper it is considered the 4-dimensional compact
homogenous space $L|_\Gamma$, where $L=H\times S^1$ is a connected Lie group,
$H$ is the Heizenberg group, $S^1$ is the circle and $\Gamma$ is the discrete
subgroup of $L$ consisting of all matrices of integer entries. The manifold $M$
generated by $L$ we equipped with an $(H,G)$-structure. Then it is a
$\W(J_1)$-manifold but it does not belong to $\W$.

In \cite{Ma} the sequel of examples begins at Examples 1 and 2 about two Engel
manifolds with almost $(H,G)$-structure. The first one is equipped with double
isotopic hyper-K\"ahlerian structures which are neither hypercomplex nor
symplectic. The second one -- with double isotopic hyper-K\"ahlerian structures
which are non-integrable but symplectic.

The Example 3 is constructed as a real semi-space with almost $(H,G)$-structure
of the class $\W$.

A real quarter-space with almost $(H,G)$-structure is shown in Example 4, which
is a $\mathcal{K}(J_1)$-manifold and an isotropic hyper-K\"ahler manifold.

In Example 5, an almost $(H,G)$-structure is introduced on a real
pseudo-hyper-cylinder in a pseudo-Euclidean real space $\mathbb{R}^5_2$. The
received manifold is not integrable and the Lie forms are non-zero, regarding
any $J_\alpha$.

The following three examples (Examples 6-8) concern several surfaces
$S^2_\mathbb{C}$ in a 3-dimensional complex Euclidean space
$\left(\mathbb{C}^3,\langle\cdot,\cdot\rangle\right)$. There is used that the
natural decomplexification of an $n$-dimensional complex Euclidean space is the
$2n$-dimensional real space with a complex skew-Hermit\-ian structure. Namely,
Example 6 shows a flat pseudo-hyper-K\"ahler manifold as a complex cylinder.
Example 7 is a complex cone with almost $(H,G)$-structure and it is a flat
hypercomplex manifold which is K\"ahlerian with respect to $J_1$ but it does
not belong to $\mathcal{W}(J_2)$ or $\mathcal{W}(J_3)$ and the Lie forms
$\theta_2$ and $\theta_3$ are non-zero. In Example 8, a complex sphere with
almost $(H,G)$-structure is given. There is shown that it is a
$\mathcal{K}(J_2)$-manifold of pointwise constant totally real sectional
curvatures, but with respect to $J_1$ and $J_3$, the corresponding Nijenhuis
tensors and Lie forms are non-zero.

The last two examples in \cite{Ma} are inspired from an example of a locally
flat almost Hermitian surface constructed in \cite{TrVa} about a connected Lie
subgroup of $\mathcal{GL}(4,\mathbb{R})$. There are introduced the appropriate
metric in two variants and then the received $(H,G)$-manifold in Example 9 is
complex with respect to $J_2$ but non-hypercomplex and the Lie forms do not
vanish; and the constructed $(H,G)$-manifold in Example 10 is flat and it is
K\"ahlerian with respect to $J_1$ but with respect to $J_2$ and $J_3$ it is not
complex.

We can use the given examples in \cite{GrMaDi1} and \cite{Ma} to construct
8-dimensional tangent bundle manifolds with almost hypercomplex
pseudo-Hermitian structure which is almost K\"ahlerian with respect to $J_1$.
For this purpose we use $J_2$ or $J_3$ like a base almost complex structure $J$
of $(M,J,g,\tilde{g})$. Then applying the above mentioned results to the case
$n=2$, we receive the following
\begin{cor} Let $M$ and $TM$ denote for short the base 4-dimensional manifold $(M,J,g,\tilde{g})$ with
almost complex skew-Hermitian structure and its tangent bundle $(TM,H,\hat{G})$
with almost hypercomplex pseudo-Hermitian structure constructed as in
\secref{sec3}, respectively.
\begin{enumerate}
    \item If $M$ with $J:=J_2$ (or $J_3$) is the complex cylinder
        introduced in Example~6 of~\cite{Ma}, then
        $TM$ is a flat pseudo-hyper-K\"ahlerian.
    \item If $M$ with $J:=J_2$ (or
        $J_3$) is the complex cone introduced in Example~7 of~\cite{Ma} or the Lie
        group introduced in Example~10 of~\cite{Ma}, then $TM$ is a flat
        almost hypercomplex pseudo-Hermitian and an
        $\AK(J_1)$-manifold.
    \item If $M$ with $J:=J_2$ is the
        complex sphere introduced in Example~8 of~\cite{Ma}, then $TM$ is
        non-flat and its $\theta_3$ is also zero as $\theta_1$.
    \item If $M$ with  $J:=J_2$ (or
        $J_3$), is some of the manifolds introduced in the other 8 known
        examples of~\cite{GrMaDi1,Ma}, then $TM$ is a non-flat
        manifold which is almost K\"ahlerian with respect to $J_1$, but
        with respect to $J_2$ and $J_3$ it is not complex and the Lie
        forms are not zero.
\end{enumerate}
\end{cor}

%%%%%%%%%%%%%%%%%%%%%%%%%%%%%%%%%%%%%%%%%%%%%%%%%%%%%%%%%%%%%%%%%%%%%%%%%%
\subsection*{Acknowledgement}
The results of this paper were obtained during the author's fellowship studies
%research stay
at Niigata University financed by Matsumae International Foundation, Japan. The
author would like to express deep gratitude to his host scientist Kouei
Sekigawa whose guidance and support were crucial for the successful completion
of this project. The author dedicate this work to the 60th anniversary of
Prof.~Sekigawa.

%%%%%%%%%%%%%%%%%%%%%%%%%%%%%%%%%%%%%%%%%%%%%%%%%%%%%%%%%%%%%%%%%%%%%%%%%%%%%%

\end{document}